\documentclass[twocolumn,noshowpacs,preprintnumbers,amsmath,amssymb,aps,eqsecnum]{revtex4}

\usepackage{graphicx}
\usepackage{dcolumn}
\usepackage{bm}

\begin{document}


\author{Marissa L. Weichman}

\affiliation{Bedford High School, 9 Mudge Way, Bedford, MA 01730}

\title{The Forgotten Night:  The Number Devil Explores Spherical Geometry}

\date{\today}

\begin{abstract}

This is a missing chapter from Hans Magnus Enzensberger's
mathematical adventure \emph{The Number Devil} (Henry Holt and
Company, New York, 1997). In the book, a math-hating boy named
Robert is visited in his dreams by the clever Number Devil, who
teaches him to love all things numerical. However, we all forget our
dreams from time to time. Here is one adventure that Enzensberger
overlooked, where the Number Devil introduces Robert to geometry
not-of-Euclid, great circles, parallel transport, the pendulum of
Foucault, and the genius of Euler.

\end{abstract}

\maketitle

The next night, Robert found himself sitting in a strange square
room with silver walls covered with rows upon rows of buttons and
switches and little blinking lights. Across the room, the Number
Devil sat perched upon a tall stool. He was peering closely at a
column of dials set into the wall. Robert cleared his throat, and
the Number Devil swiveled around on his
stool.

\bigskip

\centerline{\includegraphics[bb = 0 0 607 425,width=3.3in]{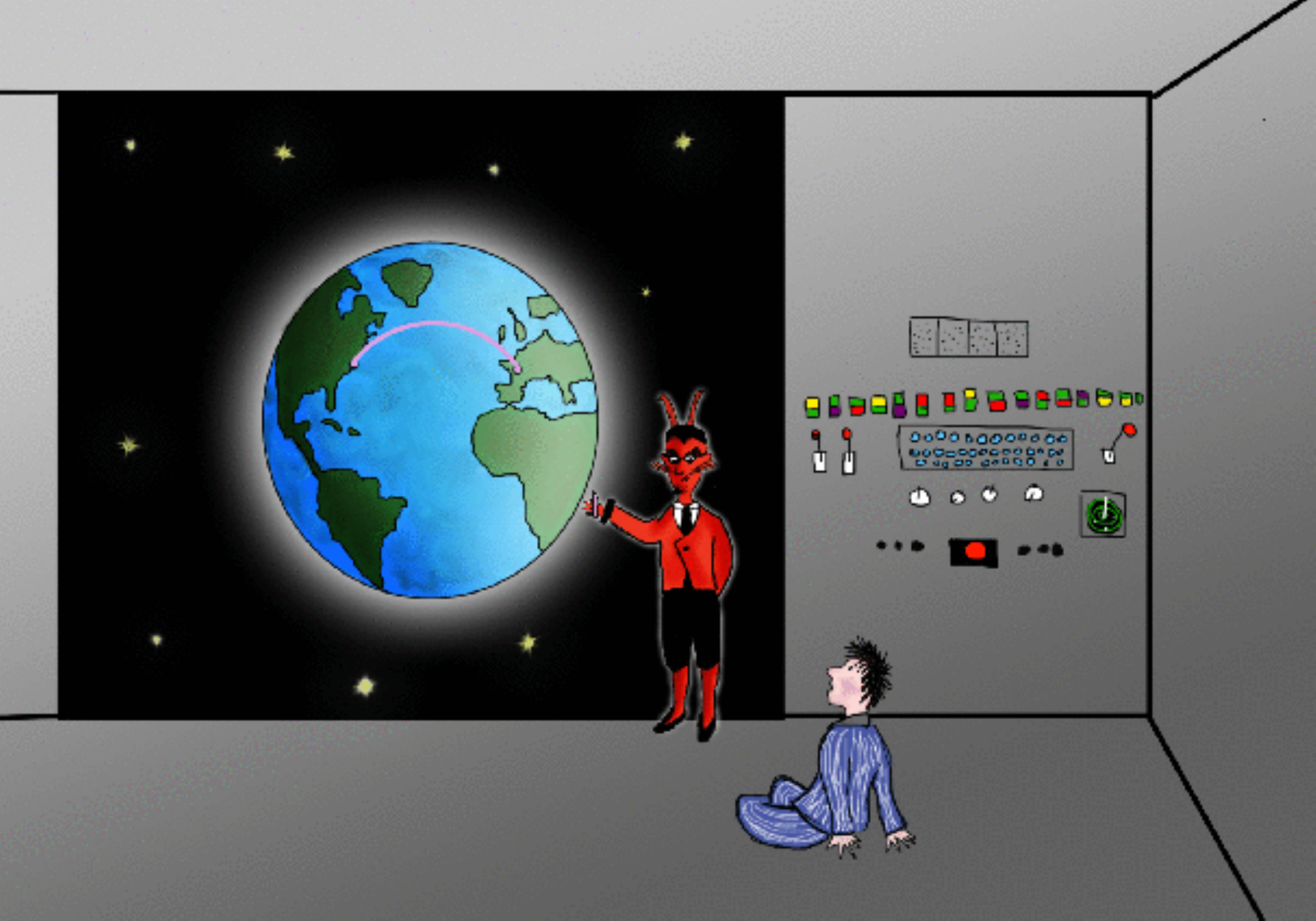}}

\bigskip

``Hello, Robert,'' he said. ``I have something quite extraordinary
to show you today. How much do you know about geometry?''

``Oh, I don't like geometry at all,'' sighed Robert. ``Geometry is
all pictures of circles and squares and parallel lines. You find how
long something is or the area of a figure, and that's it. Nothing
interesting about it.''

``How wrong you are,'' countered the Number Devil, with a
distasteful twitch of his long black whiskers. ``I'll bet in school
they only teach you Euclidean geometry.''

``What's Euclidean geometry?'' Robert asked.

``The earliest and simplest form of the geometry of planes and
three-dimensional space. The Greek mathematician who derived it,
Euclid, was a pretty clever chap. He was able to come up with all
kinds of theorems and proofs from just five simple, improvable
rules. But Euclidean geometry is pretty basic stuff. I'm going to
show you something much more interesting than just shapes on a piece
of paper.''

``But what other kind of geometry can there be?'' demanded Robert.

``Why, non-Euclidean geometry, of course!'' exclaimed the Number
Devil, jumping off his stool. ``Allow me to demonstrate.''

He pushed a large greenish-blue button on the wall, and a panel slid
back to reveal a blank whiteboard. Taking a piece of purple chalk
from his pocket, the Number Devil began to draw furiously.

\centerline{\includegraphics[bb = 0 0 378 412,width=2.7in]{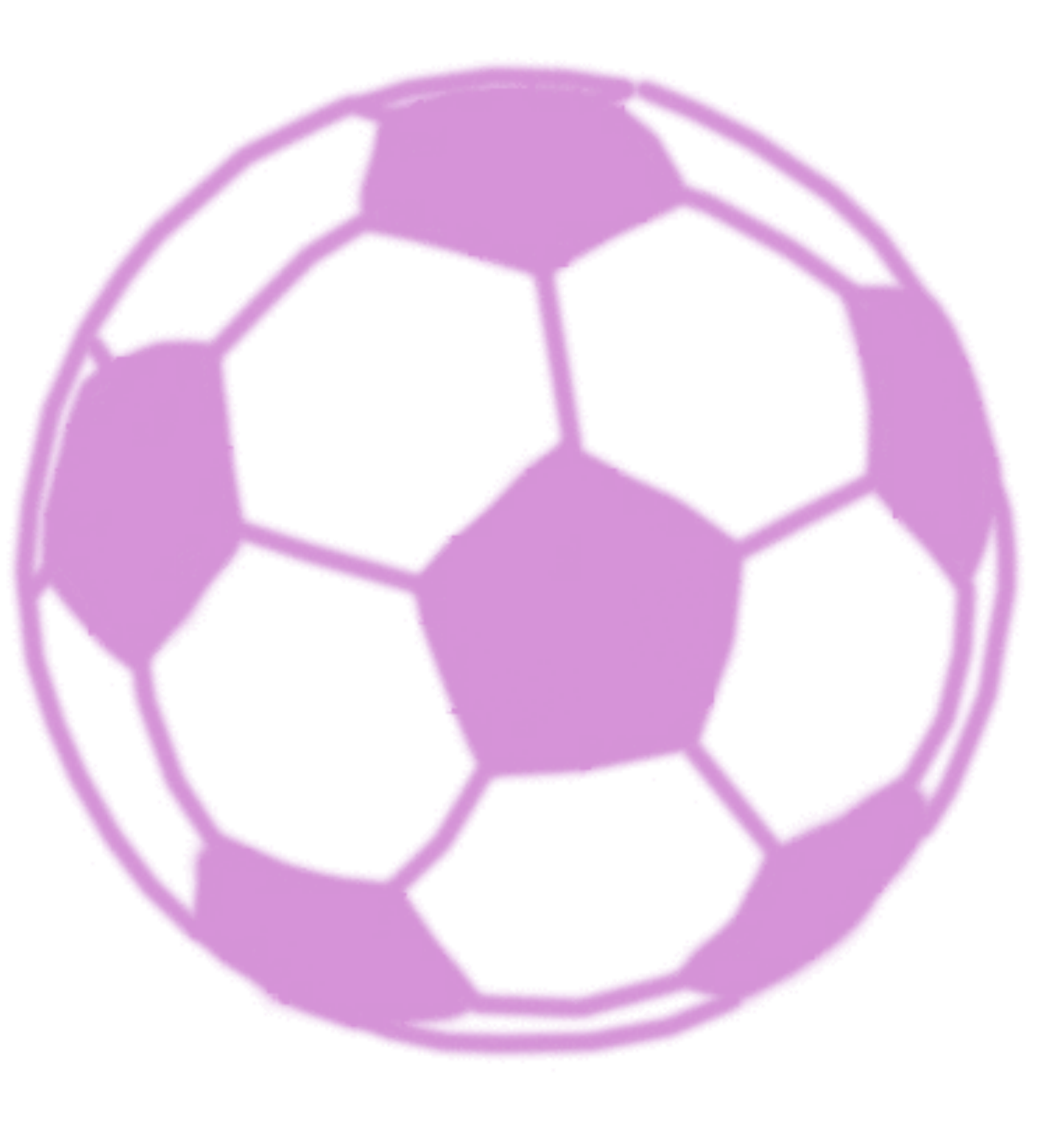}}

``Here is a soccer ball," he said, "something you probably see every
day.''

``Yes,'' said Robert.

``Can you tell me what shapes it is made of?'' asked the Number
Devil.

``It looks like pentagons and hexagons,'' said Robert. ``The purple
shapes have five sides, the white ones have six.''

``Alright,'' said the Number Devil, ``but look what happens when we
try to flatten the soccer ball out.'' With a wave of his hand, the
drawing of the soccer ball disappeared. He began scribbling
furiously again.

\centerline{\includegraphics[bb = 0 0 936 921,width=3.0in]{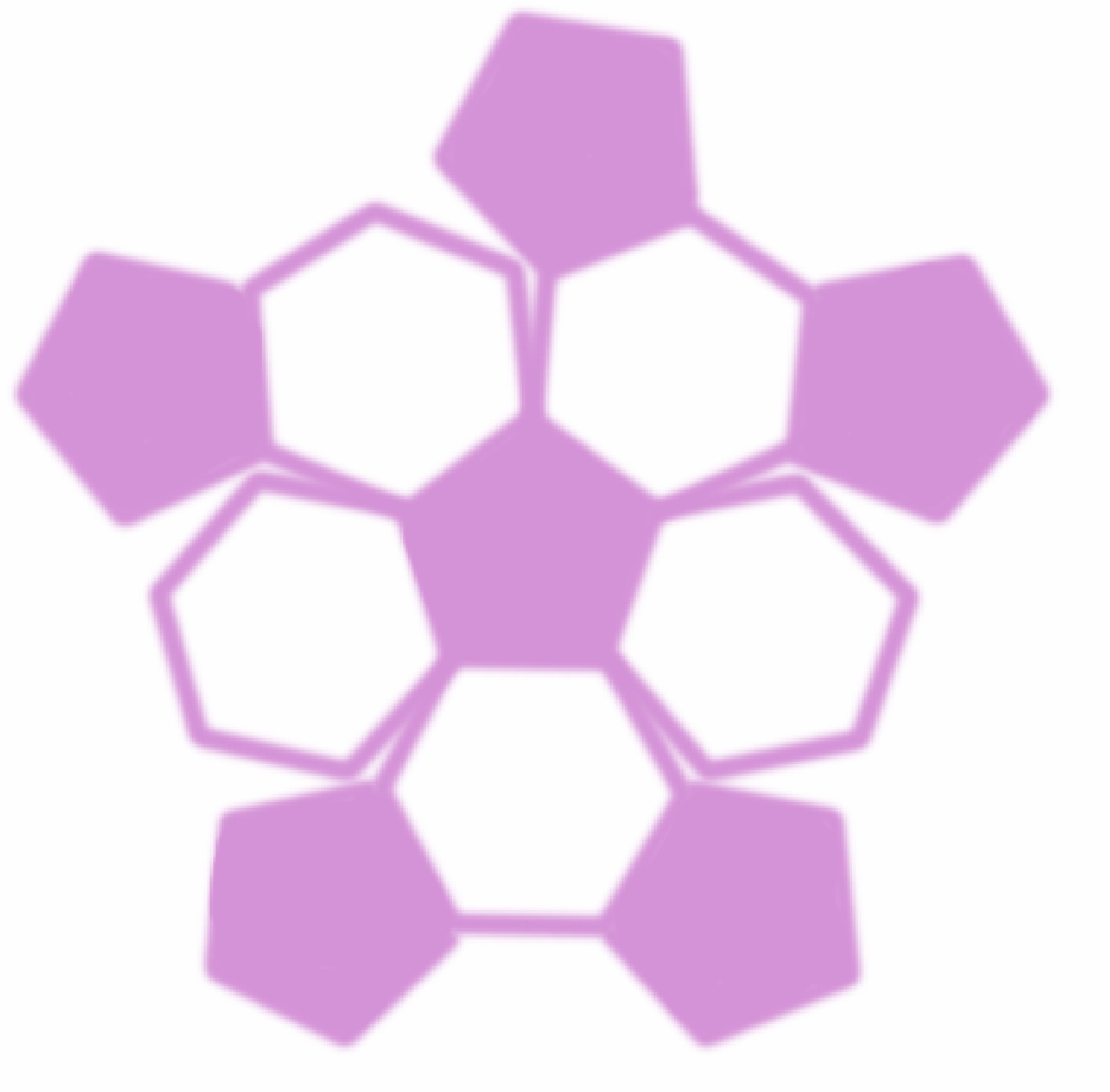}}

``On a flat surface,'' said the Number Devil, ``suddenly the
pentagons and the hexagons don't fit together anymore.''

``That can't be right,'' said Robert. ``Those aren't the same
shapes.''

``The very same,'' said the Number Devil.

``What's the catch then?''

``The catch,'' said the Number Devil, ``is that on a soccer ball we
are not dealing with Euclidean geometry. A soccer ball is a sphere,
and shapes on a sphere behave very differently than shapes on a flat
surface.''

``I think I see,'' said Robert. ``On a flat board, single interior
angles from two hexagons and one pentagon don't add up to a full 360
degrees. But they might on a sphere.''

``Exactly,'' said the Number Devil.

``That seems easy enough,'' said Robert. ``What else is different
from Euclidean geometry when you're on a sphere?''

``I thought you'd never ask!'' cried the Number Devil. He rushed
excitedly to the wall covered in buttons, where he flicked three
switches, turned a dial, and pushed a little lever up five notches.
The panels of one of the walls began to slide away. Beyond them was
a giant glass window through which Robert could see the vast reaches
of space, with little pinpoint stars scattered throughout it. The
view began to rotate slowly, and soon the giant blue-green globe of
Earth was visible.

``The Earth,'' said the Number Devil grandly. ``Quite an important
sphere to you and me, I should think.''

``Yes, I suppose,'' said Robert. ``But I don't really think of it as
a sphere when I'm living on it.''

``That's because you're so very small compared to the Earth,'' said
the Number Devil. ``At very small distances along a sphere's
surface, Euclidean geometry makes a pretty good estimate for what is
going on. But say, for example, you were in a plane traveling from
New York to Paris''---here the Number Devil marked where the two
cities were in chalk on the window---``oddly enough, you would fly
north over Newfoundland to get from point A to point B.'' He drew an
upwards-curving dashed line connecting the two purple dots.

``But why would you do that?'' Robert asked. ``It looks like you're
going far out of the way. Why wouldn't you just go across the
Atlantic in a straight line?''

``Aha!'' cried the Number Devil. ``Therein lies the problem! There
is no such thing as a straight line on a sphere. The closest thing
to a straight line is a part of a Great Circle; that is, a circle,
like the equator, that is centered on the Earth's center. And, just
like a straight line is the shortest distance between two points on
a plane, a Great Circle segment is the shortest distance between two
points on a sphere.''

``I think I understand,'' said Robert. ``If you were to take the
equator and move it so that it went through New York and Paris, it
would curve north up over Newfoundland.''

``You've got it!'' said the Number Devil.

``That's not so complicated,'' said Robert.

``But there's so much more!'' said the Number Devil. ``Say for
example, you wanted to calculate the distance between the two dots
at New York and Paris.''

``Well, if you know the radius of the Earth, you can find the
distance around a Great Circle,'' said Robert. ``That will always be
the same, won't it? It's just the circumference of a circle with the
same radius as the Earth.''

``Yes, that's right,'' said the Number Devil. ``But the problem is
finding how much of a Great Circle is being covered between those
points.''

``Okay,'' said Robert. ``So you find the angle of that slice of the
circle.''

``Yes,'' said the Number Devil. ``And there is a relatively simple
way of finding it using the latitude and longitudes of the two
cities, some trigonometry, and a lovely little device called the
Polka Dot Product.''

``Alright,'' said Robert. ``Let's see it!''

The Number Devil turned back to his whiteboard and wrote the
following:

\begin{center}
\begin{tabular}{lcc}

\emph{City} & \emph{Latitude} & \emph{Longitude} \\ \hline
\emph{NYC} & $41^\circ$ N & $74^\circ$ W \\
\emph{Paris} & $49^\circ$ N & $2^\circ$ E

\end{tabular}
\end{center}

``Now, we can use these numbers to calculate the coordinates of each
city if it were located in x-y-z space.''

``Like the unit circle, except now in three dimensions!'' said
Robert.

``Exactly,'' said the Number Devil. He erased the coordinates and
traced a large circle onto the board with two intersecting lines.
``The line running North to South is the z-axis,'' he explained.
``It runs through that black dot---which is the very center of the
Earth.''

\centerline{\includegraphics[bb = 0 0 484 534,width=2.7in]{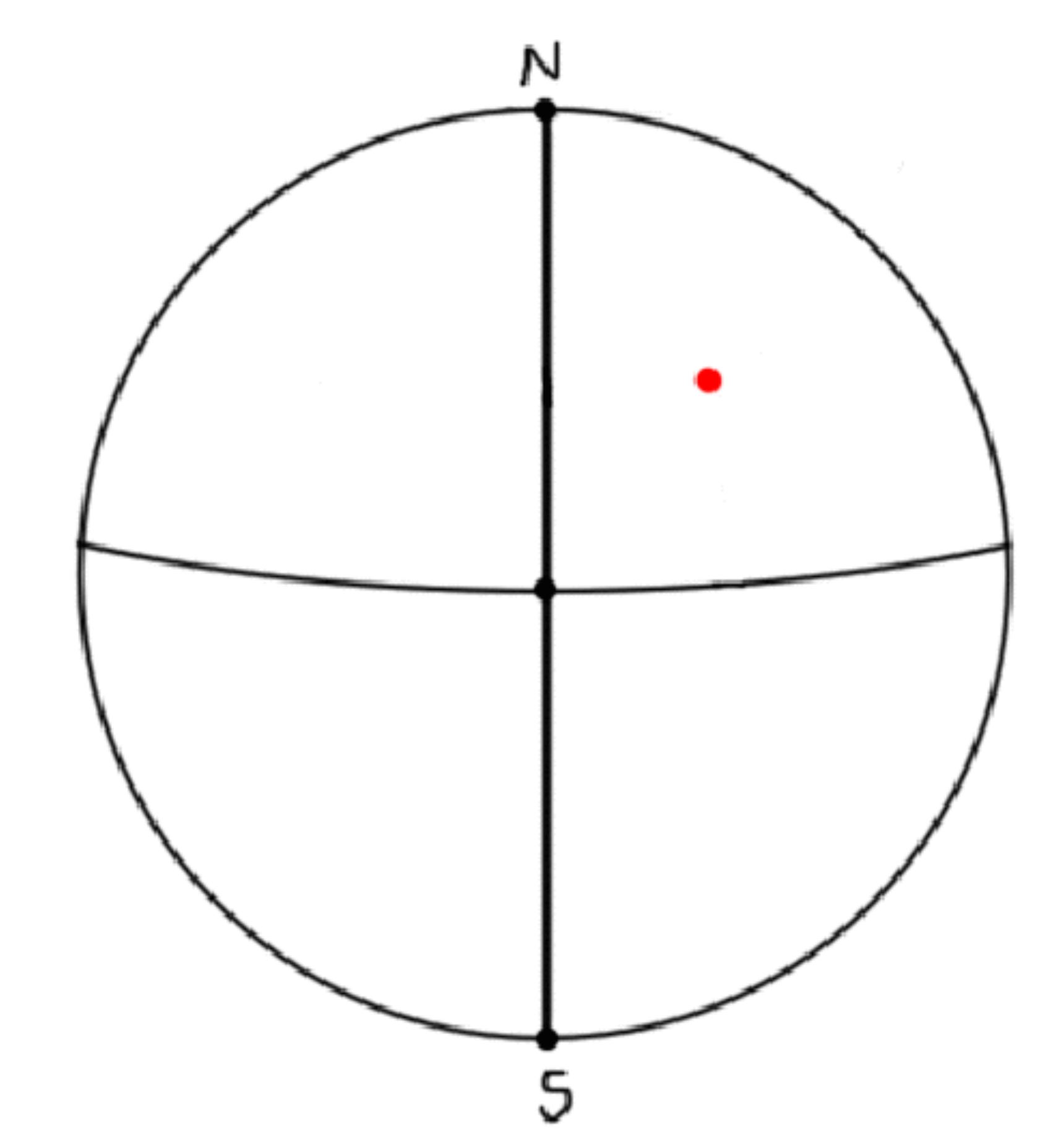}}

``Okay,'' said Robert, ``and that line from side-to-side would be
the equator.''

``Exactly,'' said the Number Devil. He then took a red piece of
chalk out of his pocket and handed it to Robert. ``Mark a point on
the circle, anywhere you like.''

``Okay,'' said Robert, ``there you go.''

``Now,'' said the Number Devil, ``that point is intersected by two
lines---one of longitude and one of latitude. Do you know the
definition of latitude, Robert?''

\centerline{\includegraphics[bb = 0 0 510 533,width=2.7in]{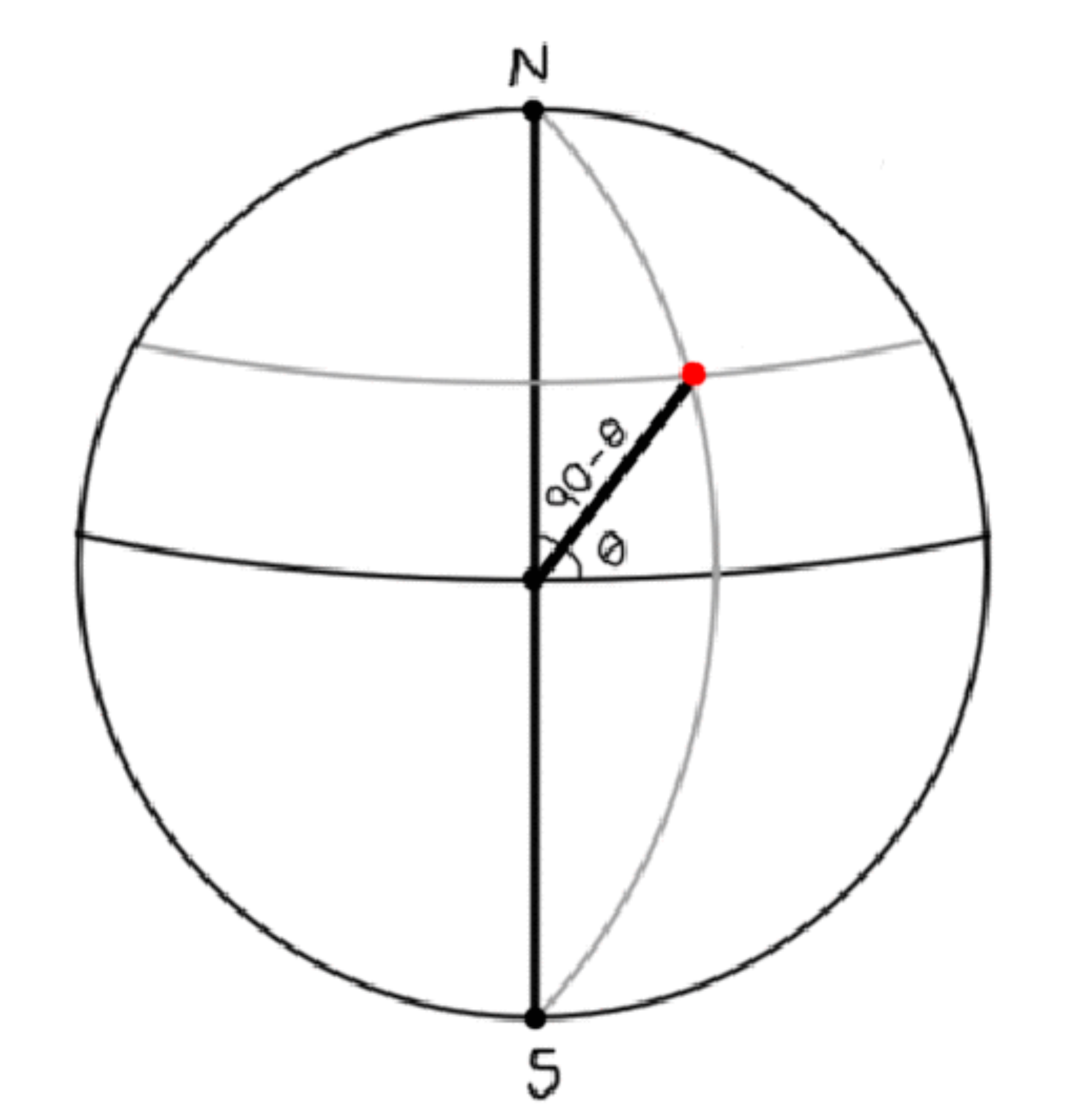}}

``Well, it must be the measure of an angle somewhere,'' said Robert.
``It's in degrees.''

``Very good,'' said the Number Devil. ``Latitude is the angle
between the equator and the line connecting your point to the center
of the Earth. Let's call it theta.''

``That would make the angle between that line and the North Pole 90
minus theta degrees'' said Robert.

``Perfect!'' declared the Number Devil.

``Okay,'' said Robert. ``We also know that the length of that line
between the red dot and the center of the Earth is just the radius
of the Earth. So if we draw in a triangle there, we could solve for
some of the coordinates.''

``You couldn't be more correct,'' said the Number Devil. ``Why don't
you go fill in as much as you can?''

``No problem,'' said Robert.

Soon the diagram looked like this:

\centerline{\includegraphics[bb = 0 0 795 417,width=3.3in]{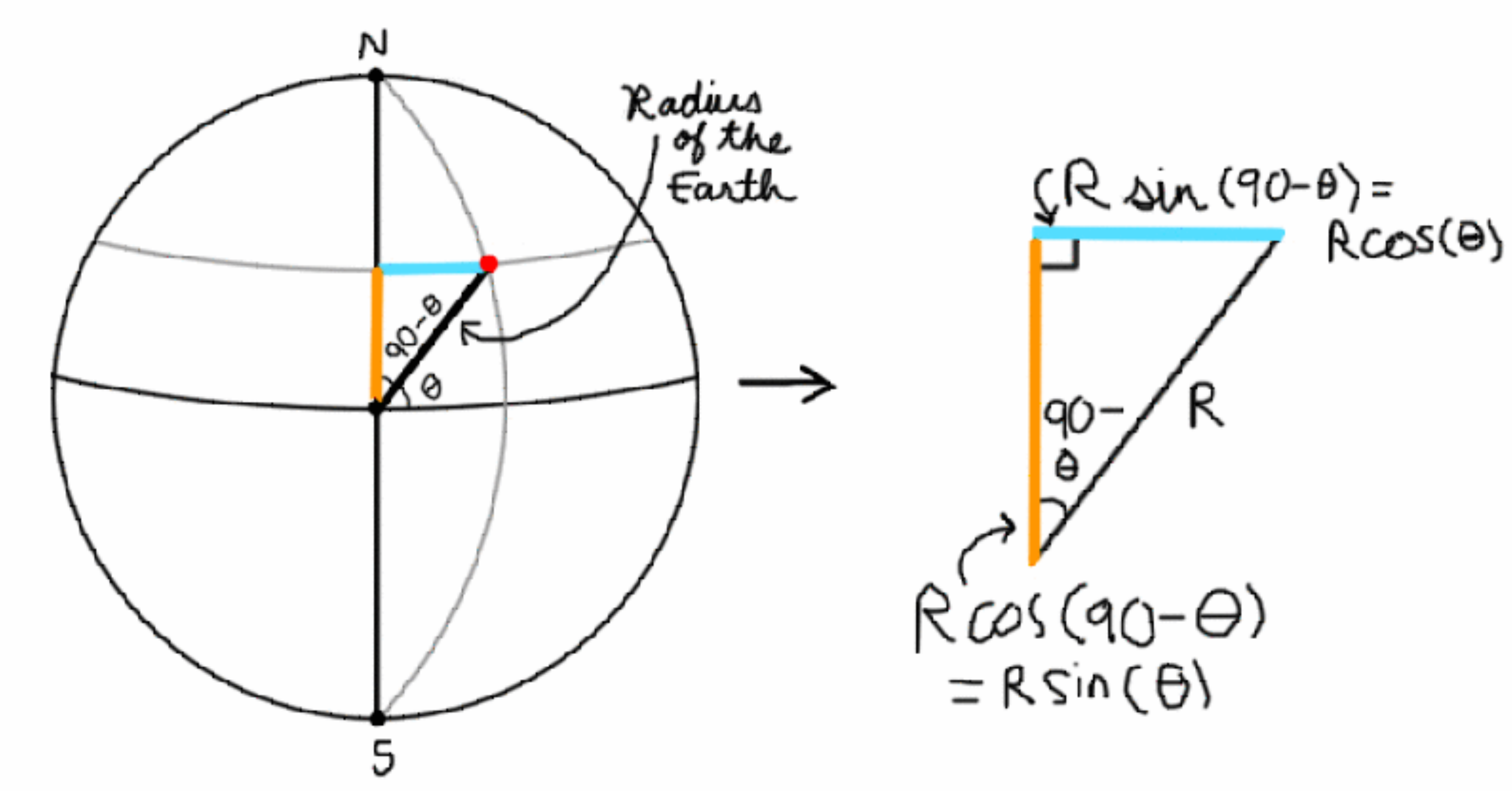}}

``Very good,'' said the Number Devil. ``I'll make a mathematician
out of you yet. Now, what do you think that orange segment can do
for us?''

``That's the height of the point!'' cried Robert. ``We've found one
of the dimensions!''

``That's right,'' said the Number Devil. ``The z-coordinate of that
$(x, y, z)$ point is the sine of the latitude times the radius of
the sphere. But in order to find the other coordinates, we'll have
to look at this slightly differently.'' The Number Devil erased the
whiteboard again, and drew another circle. ``This time, we're
looking down from the North Pole,'' he said.

\centerline{\includegraphics[bb = 0 0 677 679,width=2.8in]{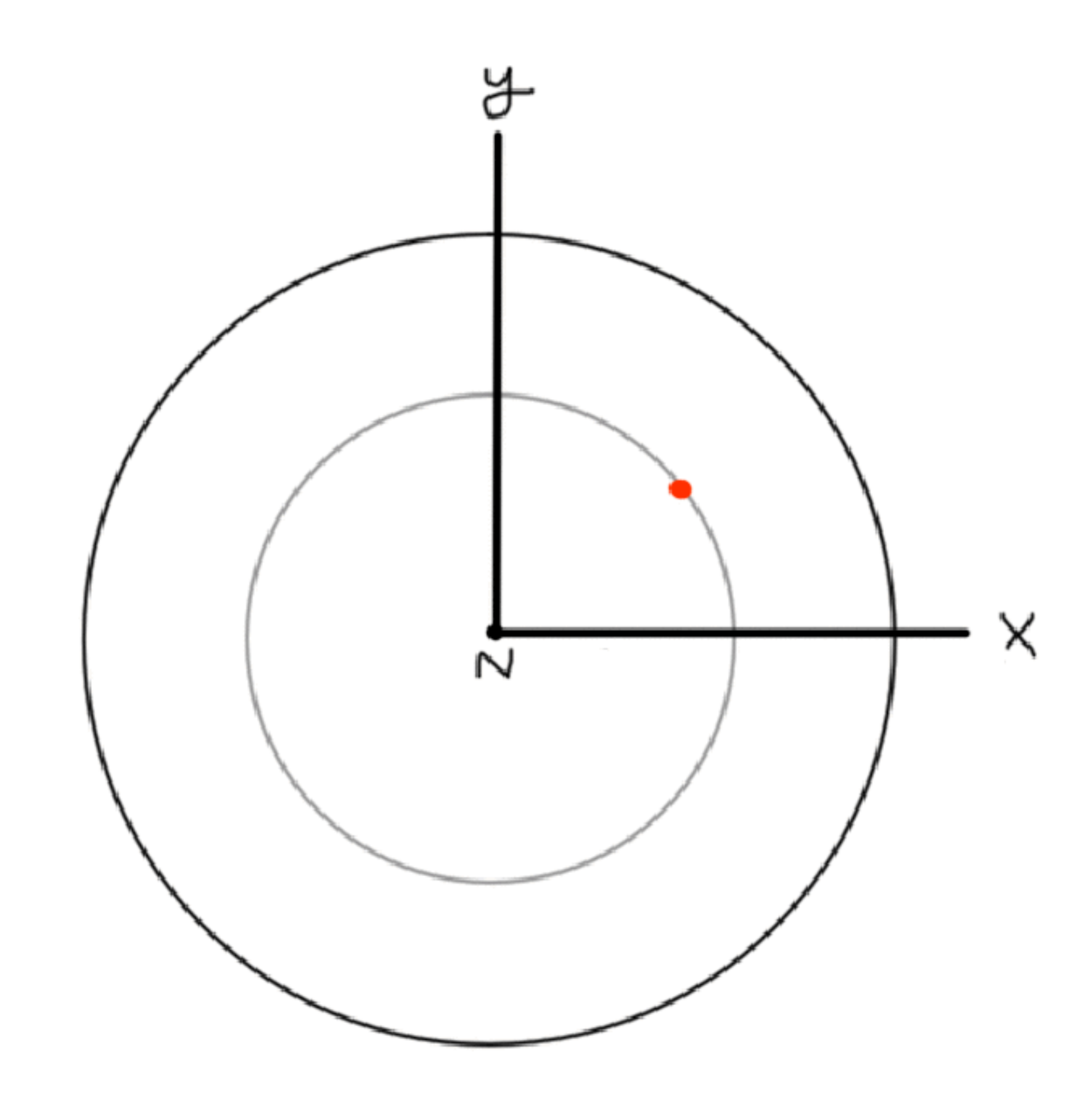}}

``I see,'' said Robert. ``The circle around that latitude is
concentric with the equator. And the x-axis looks like its covering
what would be the Prime Meridian.''

``That's right!'' said the Number Devil. ``Now, let's put our friend
longitude to use. Longitude is degrees from the prime
meridian---let's call it phi. So, draw in a triangle, Robert.''

\centerline{\includegraphics[bb = 0 0 510 512,width=2.8in]{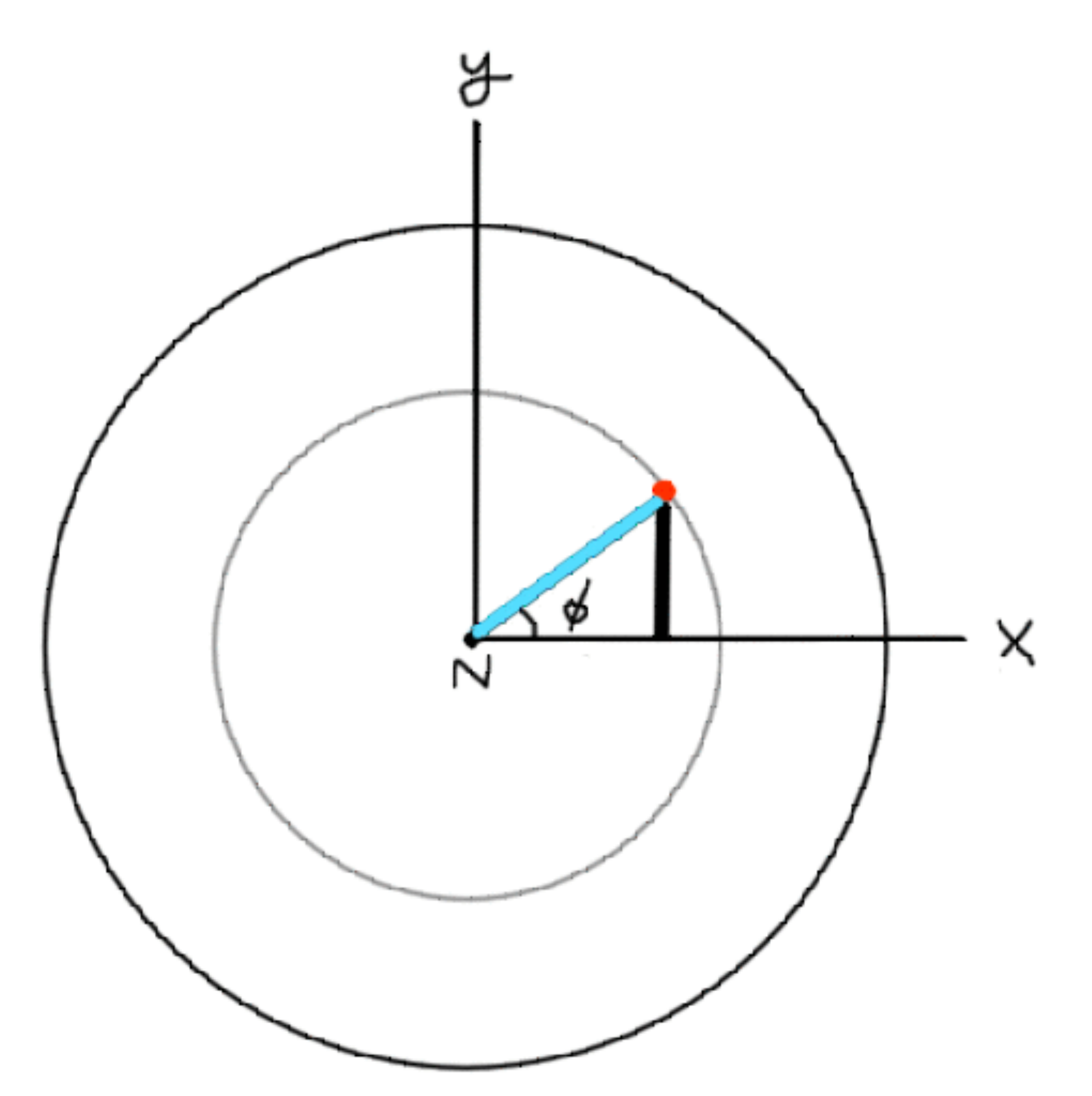}}

``Wait half a second,'' said Robert. ``This time the radius of the
Earth doesn't come into the picture. How are we supposed to find the
lengths of the dimensions?''

``I'm so glad you asked,'' said the Number Devil. ``Let's think back
to our first circle. Remember that pale blue line that connected the
red point directly to the z-axis, making a right angle?''

``I think I see the connection,'' said Robert. ``If we add the
triangle to this diagram, the point connecting the North-South axis
to the red dot is that same length, the radius of the Earth times
the cosine of theta!''

``Brilliant,'' said the Number Devil. ``Now in fill everything.''

\centerline{\includegraphics[bb = 0 0 507 327,width=2.7in]{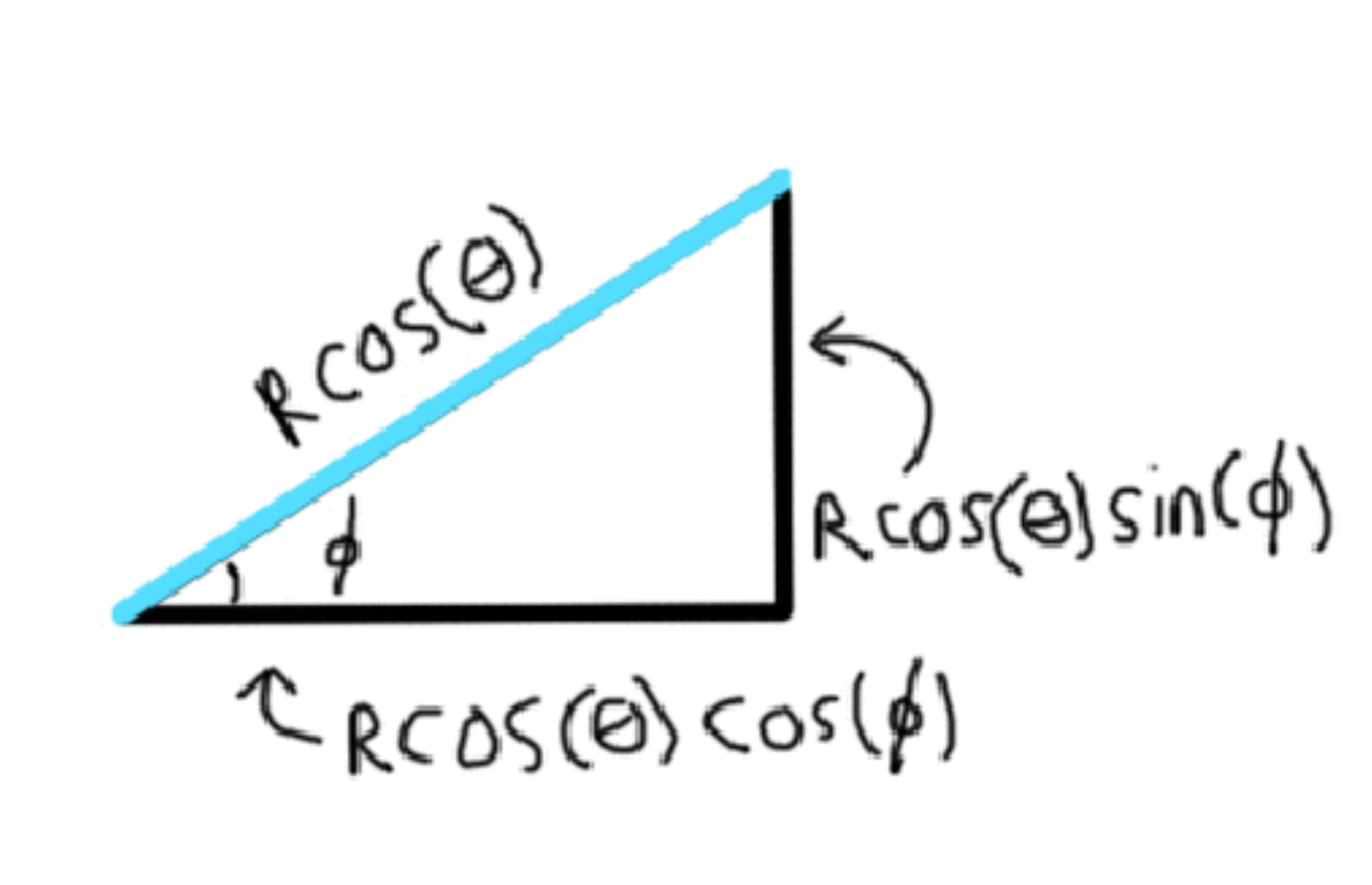}}

``There,'' said Robert. ``And that gives us the x- and y-coordinates
too! The altitude of that triangle is along the y-axis and the base
is along the x-axis. So here is what we know\ldots''

\bigskip

\begin{center}
\begin{minipage}{2.9in}
\emph{A point with latitude and longitude of $(\theta,\phi)$ has
three-dimensional coordinates:}

\smallskip

\noindent $[R\cos(\theta) \cos(\phi), R\cos(\theta)\sin(\phi),
R\sin(\theta)]$
\end{minipage}
\end{center}

\bigskip

``Lovely,'' said the Number Devil. ``Now, given that the radius of
the Earth is around 6378 kilometers, finding the coordinates of New
York and Paris is just a matter of plugging-and-chugging.''

``I can handle that,'' said Robert. He wrote out:

\begin{eqnarray}
\underline{NYC:} && \theta = 49^\circ,\ \ \phi = -74^\circ
\nonumber \\
x &=& 6378 \cos(41^\circ) \cos(-74^\circ) = 1327\ \mathrm{km}
\nonumber \\
y &=& 6378 \cos(41^\circ) \sin(-74^\circ) = -4627\ \mathrm{km}
\nonumber \\
z &=& 6378 \sin(41^\circ) = 4184\ \mathrm{km}
\nonumber \\
&&\to (1327, -4627, 4184)
\nonumber
\end{eqnarray}

``Now for Paris,'' said Robert.

\begin{eqnarray}
\underline{Paris:} && \theta = 49^\circ,\ \ \phi = 3^\circ
\nonumber \\
x &=& 6378 \cos(49^\circ) \cos(3^\circ) = 4179\ \mathrm{km}
\nonumber \\
y &=& 6378 \cos(49^\circ) \sin(3^\circ) = 219\ \mathrm{km}
\nonumber \\
z &=& 6378 \sin(49^\circ) = 4814\ \mathrm{km}
\nonumber \\
&&\to (4179, 219, 4814)
\nonumber
\end{eqnarray}

``So we have two sets of coordinates in three dimensions,'' said
Robert, ``but we still don't know how to find the angle formed
between the center of the sphere and those two points.''

``Think about it this way,'' said the Number Devil, turning to the
whiteboard. "We now have a triangle with these three points, and we
know the coordinates of all three vertices. It's no trouble at all
to find one of the angles.''

\centerline{\includegraphics[bb = 0 0 679 475,width=2.9in]{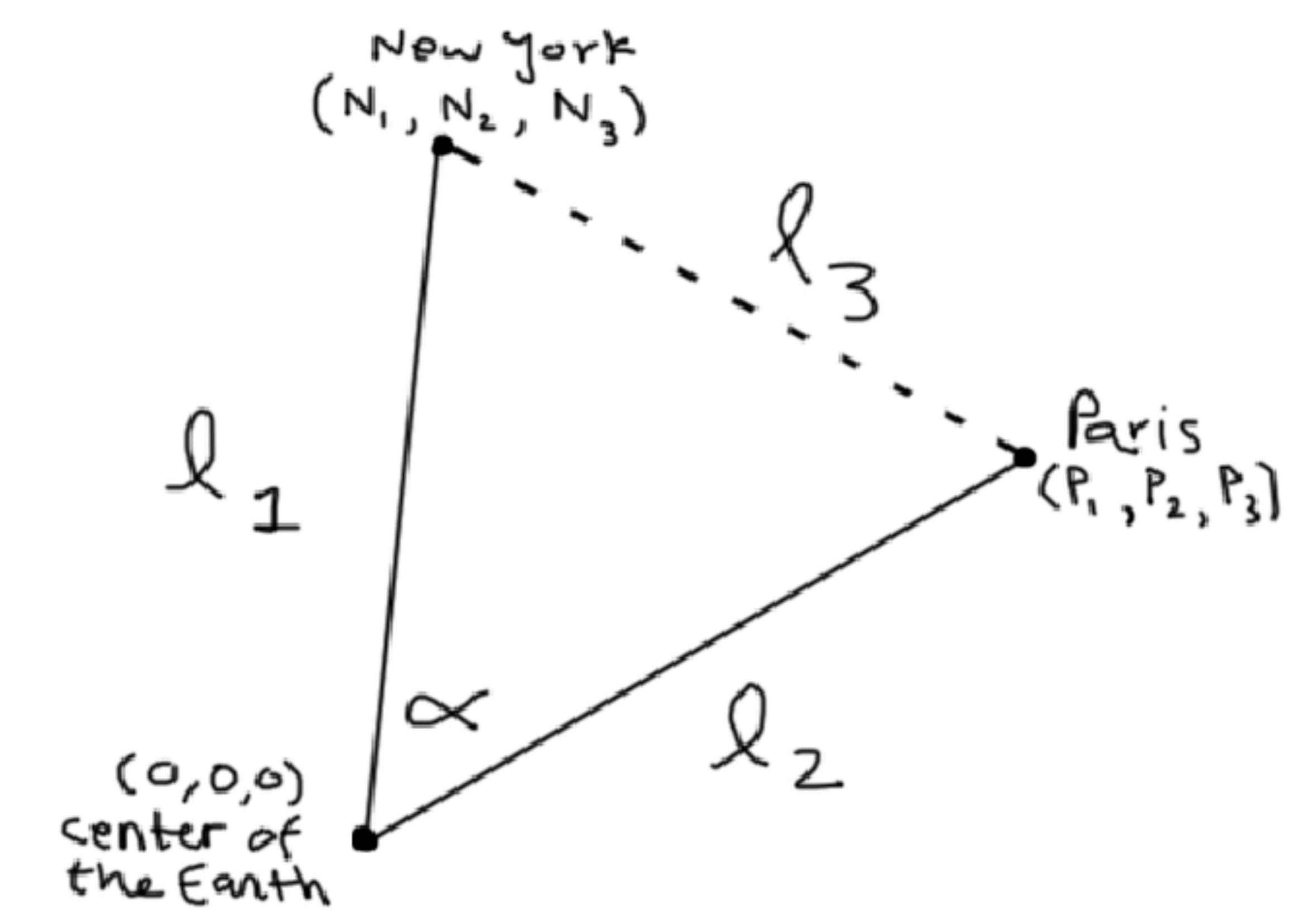}}

``I see,'' said Robert. ``You can use the distance formula to
calculate the lengths of the three sides, and then use the cosine
formula to solve for alpha, the angle we need!''

``Perfect thinking, my boy,'' said the Number Devil.

Robert grabbed a piece of chalk and wrote out:

\begin{eqnarray}
l_1 &=& \sqrt{N_1^2 + N_2^2 + N_3^2}
\nonumber \\
l_2 &=& \sqrt{P_1^2 + P_2^2 + P_3^2}
\nonumber \\
l_3 &=& \sqrt{(N_1-P_1)^2 + (N_2-P_2)^2 + (N_3-P_3)^2}
\nonumber
\end{eqnarray}

``Now,'' said the Number Devil, ``you know the cosine formula, don't
you?''

``Sure,'' said Robert.

\begin{eqnarray}
&&l_3^2 = l_1^2 + l_2^2 - 2 l_1 l_2 \cos(\alpha)
\nonumber \\
&\to& \cos(\alpha) = \frac{l_1^2 + l_2^2 - l_3^2}{2 l_1 l_2}
\nonumber
\end{eqnarray}

``Now we just need to do a bit of substitution,'' said the Number
Devil.

``Okay,'' said Robert. ``The numerator of that fraction looks like
it can be simplified quite a bit\ldots''

\begin{eqnarray}
l_1^2 + l_2^2 - l_3^2 &=& N_1^2 + N_2^2 + N_3^2 + P_1^2 + P_2^2 + P_3^2
\nonumber \\
&-& (N_1-P_1)^2 - (N_2-P_2)^2 - (N_3-P_3)^2
\nonumber \\
&=& N_1^2 + N_2^2 + N_3^2 + P_1^2 + P_2^2 + P_3^2
\hskip-1.97in \rule[0.042in]{1.97in}{0.01in}  
\nonumber \\
&&-\ (N_1^2 + N_2^2 + N_3^2 + P_1^2 + P_2^2 + P_3^2)
\hskip-2.05in\rule[0.042in]{2.05in}{0.01in}   
\nonumber \\
&&+\ 2(N_1 P_1 + N_2 P_2 + N_3 P_3)
\nonumber \\
&=& 2(N_1 P_1 + N_2 P_2 + N_3 P_3)
\nonumber
\end{eqnarray}

``That looks much prettier,'' said Robert. ``Now the cosine formula
is just:''

\begin{equation}
\cos(\alpha) = \frac{N_1 P_1 + N_2 P_2 + N_3 P_3}{l_1 l_2}
\nonumber
\end{equation}

``And that,'' said the Number Devil, ``is our friend the Polka Dot
Product. It is defined as\ldots''

\begin{eqnarray}
{\bf x}_1 = (x_1,y_1,z_1),\ \ {\bf x}_2 = (x_2,y_2,z_2)
\nonumber \\
{\bf x}_1 \cdot {\bf x}_2 = x_1 x_2 + y_1 y_2 + z_1 z_2.
\nonumber
\end{eqnarray}

``So therefore,'' the Number Devil continued, ``our cosine formula
is just:''

\begin{equation}
\cos(\alpha) = \frac{{\bf l}_1 \cdot {\bf l}_2}{l_1 l_2}
\nonumber
\end{equation}

``Is that elegant or what!'' cried the Number Devil.

``It is pretty nice,'' Robert admitted. ``Now it's no hassle at all
to plug in our coordinates for New York City and Paris! We already
know $l_1$ and $l_2$ are both just the radius of the Earth, or 6378
kilometers\ldots''

\begin{eqnarray}
NYC &=& (1327, -4627, 4184)
\nonumber \\
Paris &=& (4179, 219, 4814)
\nonumber \\
\cos(\alpha) &=& \frac{(1327 \cdot 4179) + (-4627 \cdot 219)
+ (4184 \cdot 4814)}{6378 \cdot 6378}
\nonumber \\
&=& 0.6065 \ \to \ \alpha = \cos^{-1}(0.6065) = 52.66^\circ
\nonumber
\end{eqnarray}

``Finally!'' cried Robert. ``We found the angle. Now finding the
distance between New York and Paris is as easy as pie!''

\begin{equation}
\frac{52.66}{360} \cdot 2\pi (6378) = 5862\ \mathrm{km}!!!
\nonumber
\end{equation}

``Phew!'' sighed Robert. ``That was a lot of work for one silly
little distance. I'm almost ready for a nap.''

``But there are so many more interesting things to do with
spheres!'' the Number Devil exclaimed.

``Like what?'' Robert yawned.

``Almost anything you want!'' exclaimed the Number Devil.

``What about area?'' asked Robert. ``I bet the area of my backyard
could be estimated as if it was on a flat plane, but that probably
wouldn't be accurate for points that are far apart.''

``You couldn't be closer to the truth, Robert,'' said the Number
Devil proudly. ``The area of a shape on a sphere goes back to our
soccer ball problem. Remember how the corners of the hexagons and
pentagons fit together on a sphere, but didn't when they were laid
out flat?''

``Yes,'' said Robert.

``In the same way, the three angles of a flat triangle add up to 180
degrees, but if you had a triangle on a sphere, that wouldn't
necessarily be the case.''

``That makes sense,'' said Robert. ``But what do the angles have to
do with the area of the triangle?''

``It turns out that the area of the figure is directly related to
the difference between the sum of its angles, and the sum it should
have if it were on a plane. So, for example, the area for a triangle on
a sphere would be:''

\begin{equation}
A = R^2 (\alpha + \beta + \gamma - 180^\circ) \cdot
\left(\frac{2\pi}{360} \right)
\nonumber
\end{equation}

``That seems simple enough,'' said Robert. ``But given only the
latitudes and longitudes of the corners of the triangle, how can you
find those angles?''

``That's where it gets a bit tricky," said the Number Devil. "We
have to use both the Polka Dot Product and its friend, the Criss
Cross Product.''

``I'm up for that!'' said Robert. ``I'd like to find the area of the
Bermuda Triangle.''

``Fantastic!'' said the Number Devil. ``The three corners of the
Bermuda Triangle are located in Bermuda (naturally), Melbourne,
Florida and Puerto Rico. Their latitudes and longitudes are as
follows\ldots''

\begin{center}
\begin{tabular}{lcc}
& \emph{Latitude} & \emph{Longitude} \\ \hline
\emph{Florida} & $28^\circ$ N & $81^\circ$ W \\
\emph{Bermuda} & $32^\circ$ N & $65^\circ$ W \\
\emph{Puerto Rico} & $18^\circ$ N & $66^\circ$ W
\end{tabular}
\end{center}

``Now we just use the formula we used earlier to turn those into
three-dimensional coordinates,'' said Robert. ``So they come out
like this\ldots''

\begin{eqnarray}
\mbox{\it Florida} &=& (881, -5562, 2994)
\nonumber \\
\mbox{\it Bermuda} &=& (2286, -4902, 3380)
\nonumber \\
\mbox{\it Puerto\ Rico} &=& (2467, -5541, 1971)
\nonumber
\end{eqnarray}

``Okay,'' said the Number Devil. ``That's all fine and dandy, but
before we put those to use, I'm going to need to introduce you to
our friend the Criss Cross Product. Now, the sides of a triangle on
a sphere are all segments of Great Circles---that's the definition
of the triangle. So, each corner of the triangle is the intersection
of two Great Circles.''

\centerline{\includegraphics[bb = 0 0 373 378,width=2.5in]{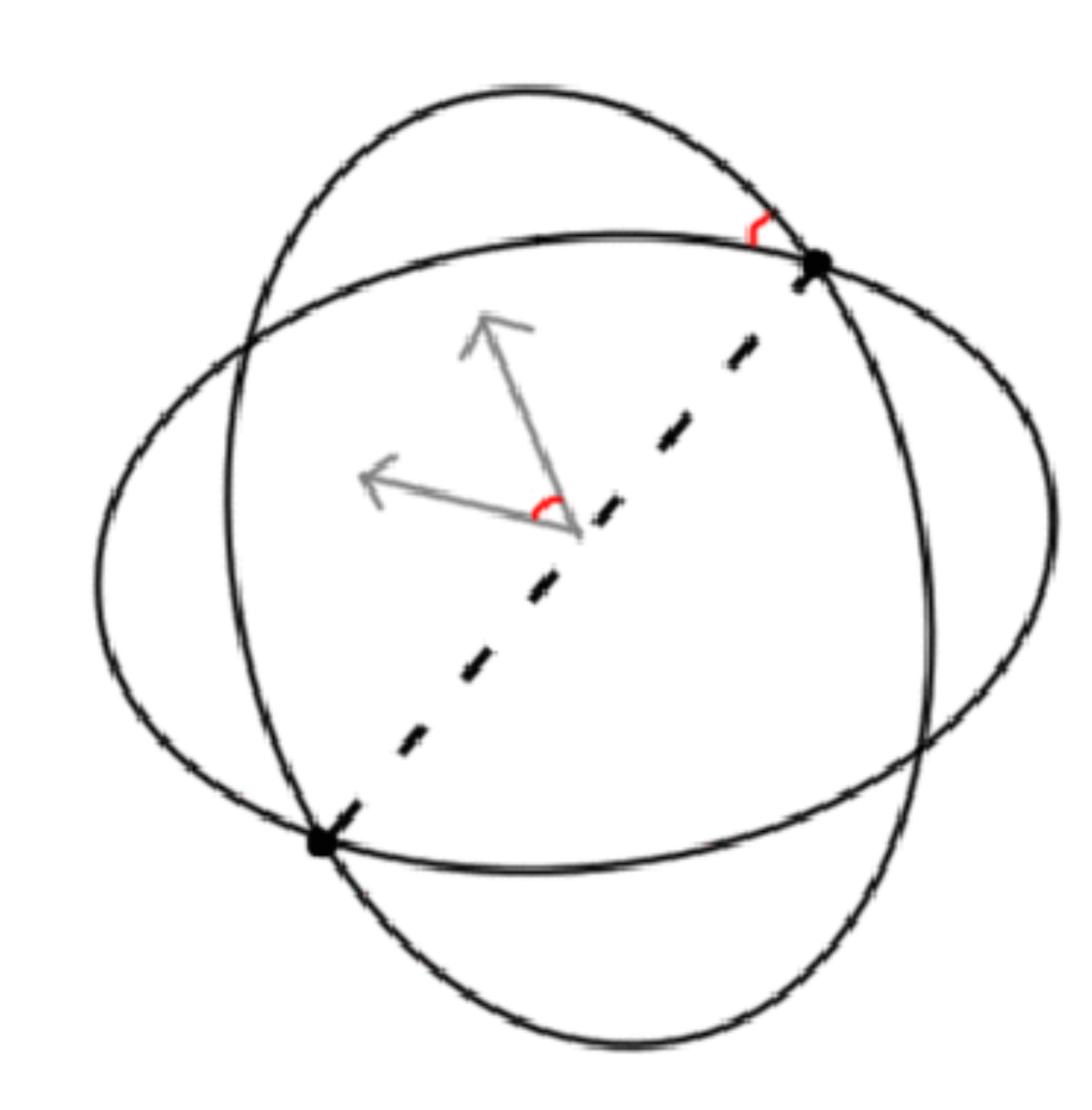}}

``And we're trying to find the angle where those two Great Circles
intersect!'' said Robert. ``That makes sense.''

``Now,'' said the Number Devil. ``This is a bit difficult to
picture, but imagine if you drew two lines perpendicular to each of
the Great Circles.''

``I think I understand,'' said Robert. ``The angle between the two
perpendiculars would be congruent to the angle between the two
circles!''

``That's correct!'' said the Number Devil.  ``So for a triangle on a
sphere, there are three Great Circles, each passing through two
corners of the triangle. Therefore, there are three perpendiculars
to be found. And our friend, the Criss Cross Product, can give us the
perpendicular to any plane!''

\centerline{\includegraphics[bb = 0 0 581 306,width=3.0in]{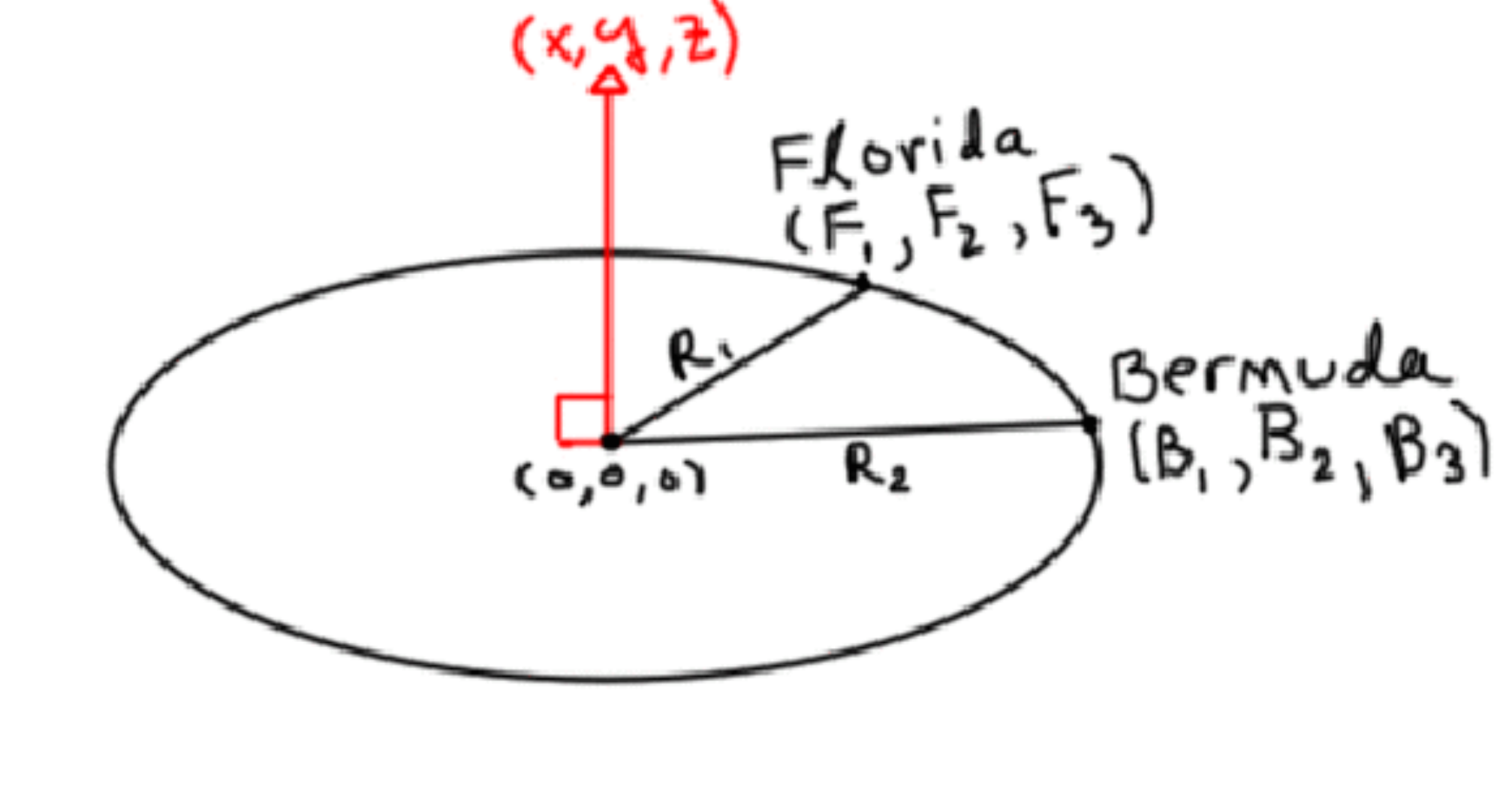}}

``And then once we have the three perpendiculars,'' said Robert,
``we can just use the Polka Dot Product to find the angle between
any two of them!''

``Spot on!'' said the Number Devil. ``Now, say we're looking at the
Great Circle containing our points in Florida and Bermuda, and our
perpendicular exits the sphere at (x, y, z). The Criss Cross Product
is like the determinant of a matrix\ldots''

\begin{eqnarray}
{\bf F} \times {\bf B}
&=& \left|\begin{array}{ccc}
x & y & z \\
F_1 & F_2 & F_3 \\
B_1 & B_2 & B_3
\end{array} \right|
\nonumber \\
&=& x(F_2 B_3 - B_2 F_3)
\nonumber \\
&&+\ y(F_3 B_1 - B_3 F_1)
\nonumber \\
&&+\ z(F_1 B_2 - B_1 F_2)
\nonumber
\end{eqnarray}

``But instead of computing the entire determinant like I have
above,'' the Number Devil continued, ``you just take the parts being
multiplied by x, y and z, and put them into one three-dimensional
coordinate!''

\begin{widetext}
\begin{equation}
\left. \begin{array}{c}
x(F_2 B_3 - B_2 F_3) \\
y(F_3 B_1 - B_3 F_1) \\
z(F_1 B_2 - B_1 F_2)
\end{array} \right\}
\to (F_2 B_3 - B_2 F_3, F_3 B_1 - B_3 F_1, F_1 B_2 - B_1 F_2)
\nonumber
\end{equation}
\end{widetext}

``And, if you scale the length to be the radius of the Earth, you
have the point where the perpendicular leaves the sphere!'' cried
the Number Devil.

``So you just calculate that point for each pair of points, and then
take the Polka Dot Products of each pair of the perpendiculars to
find each angle,'' said Robert. ``That sounds like entirely too much
algebra to me.''

``It is an awful lot of computation,'' said the Number Devil.
``Luckily, I have a nice little computer program that will do all
the math for us, and spare us the possibility of making mistakes.''
He walked over to a panel in the wall covered with glowing buttons
and punched a few rapidly. Robert heard three loud beeps,
and suddenly rows of numbers appeared on the screen.

``These tell us the longitudes and latitudes of the points where the
perpendiculars emerge from the sphere,'' the Number Devil explained.

\begin{eqnarray}
\mbox{\sl FLORIDA} \times \mbox{\sl PUERTO RICO} &\to&
(48^\circ\ \mathrm{N},45^\circ\ \mathrm{E})
\nonumber \\
\mbox{\sl PUERTO RICO} \times \mbox{\sl BERMUDA} &\to&
(3^\circ\ \mathrm{N},157^\circ\ \mathrm{W})
\nonumber \\
\mbox{\sl BERMUDA} \times \mbox{\sl FLORIDA}  &\to&
(56^\circ\ \mathrm{S},43^\circ\ \mathrm{W}).
\nonumber
\end{eqnarray}

``And now with one more push of a button,'' said the Number Devil,
``I can have the computer print out the angles given once the Polka
Dot Products are taken!''

\begin{center}
\begin{tabular}{lcc}
{\sl FLORIDA} & $\to$ &$52.8^\circ$ \\
{\sl PUERTO RICO} & $\to$ & $54.8^\circ$ \\
{\sl BERMUDA} & $\to$ & $74.1^\circ$ \\ \hline
{\sl TOTAL} & $\to$ & $181.7^\circ$
\end{tabular}
\end{center}

``And look at that!'' cried Robert. ``They add up to more than 180
degrees''

``See, I told you,'' said the Number Devil. ``Spheres are tricky
blighters. Now, Robert, calculate the area. You don't need a fancy
computer for that.''

``Okay,'' said Robert, turning back to the whiteboard. ``Shouldn't
be any problem at all.''

\begin{equation}
6378^2 \cdot \left[181.7^\circ - 180^\circ \right] \cdot
\left(\frac{2\pi}{360} \right) = 1,211,500 \mbox{ km}^2
\nonumber
\end{equation}

``Now, just for the sake of comparison,'' said the Number Devil,
``let's look at what the area of this triangle would be were it
flat''

``Alright,'' said Robert. ``We can find the Great Circle distances
between Florida, Bermuda and Puerto Rico, just like we did for New
York and Paris. That would give us\ldots''

\begin{eqnarray}
a &=& \mbox{{\it Florida -- Puerto Rico}: 1895 km}
\nonumber \\
b &=& \mbox{{\it Puerto Rico --  Bermuda}: 1562 km}
\nonumber \\
c &=& \mbox{{\it Bermuda -- Florida}: 1604 km}
\nonumber
\end{eqnarray}

``Very good,'' said the Number Devil. ``Now, I assume you know of
the Heroic Formula?''

``Sure,'' said Robert.

\begin{equation}
A = \sqrt{s(s-a)(s-b)(s-c)},\ \ s = \frac{a+b+c}{2}
\nonumber
\end{equation}

``And when we plug in our side lengths,'' Robert continued, ``we get
the area to be 1,200,800 km$^2$. That's around 10,700 km$^2$extra
real estate!''

``Unfortunately it's mostly ocean,'' said the Number Devil. ``But
see what a difference the sphere makes, even with relatively small
triangles?''

``It is pretty neat,'' Robert admitted.

``So, you wouldn't mind if I were to show you another way to
calculate the sum of the angles of a polygon on a sphere?'' asked
the Number Devil.

``Fine,'' said Robert. ``Although I am getting very sleepy.''

``Alright,'' said the Number Devil. ``Follow me, please.'' He led
Robert through a door that had suddenly appeared in the wall and
down a long silvery corridor. At the end of the hall was yet another
door, but this one opened out into the vast emptiness of space. In
front of Robert was a bridge tethered to the vessel he had just left
(which he now saw was a space station) that led to a small round
asteroid.

The Number Devil then handed Robert a strange device from one of his
pockets---it was a wheel with little notches marking degrees, like a
protractor. The device was mounted on a perfectly frictionless axle,
so even when the axle was twisted back and forth, the wheel would
remain motionless.

\centerline{\includegraphics[bb = 0 0 191 229,width=2.5in]{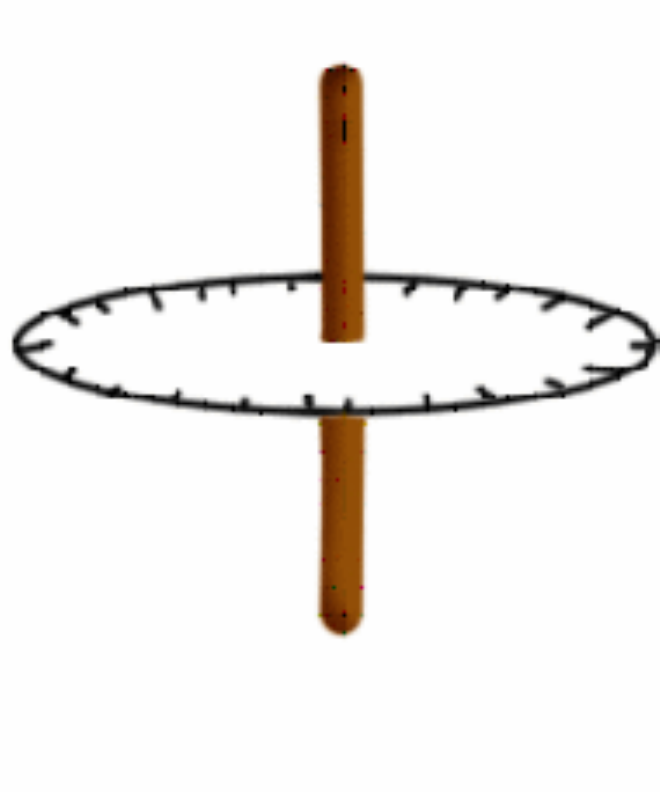}}

``This is called a parallel transporter!'' said the Number Devil.
``Go stand at the very top of that asteroid, keeping the wheel held
parallel with the ground.''

Robert obliged.

``Now, see where the zero degree marking is on the wheel?'' asked the
Number Devil. ``Walk in that direction for three steps.''

\centerline{\includegraphics[bb = 0 0 355 540,width=3.0in]{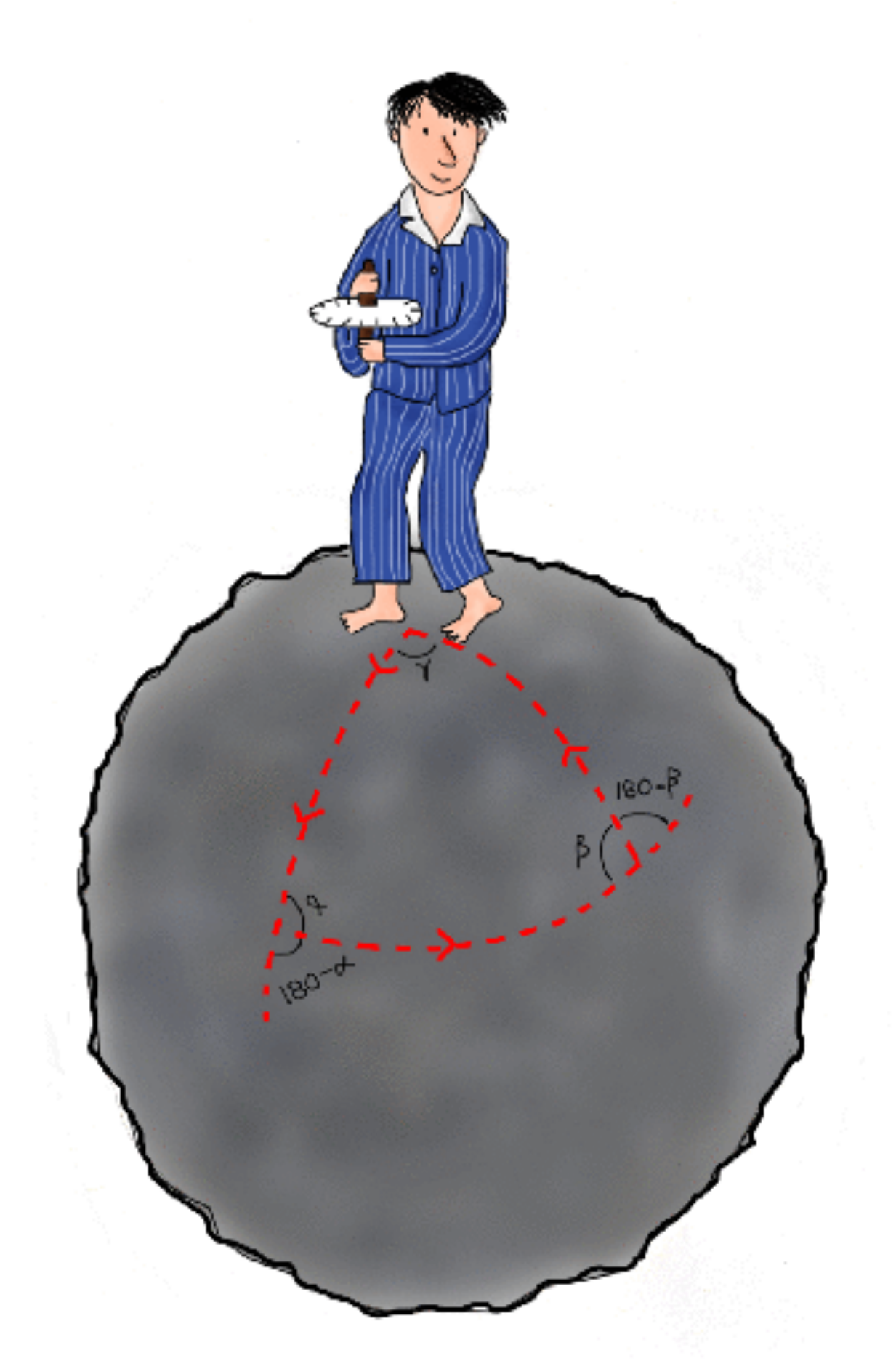}}

Again, Robert did as told.  He was now almost standing sideways on
the asteroid!

``Make a left turn. The zero degree marker should still be facing
the direction you were walking. Say the angle you turned was of
alpha degrees. Which marker on the wheel is lined up with your new
direction of travel?''

``That would be 180 minus alpha!'' said Robert.

``That's right,'' said the Number Devil. ``Now take three more steps
and make another turn, this time of angle beta.''

Robert had to look up to see the Number Devil, since his head was
pointed almost toward him.  ``Now the wheel is reading 180 minus
alpha, plus 180 minus beta!'' said Robert.

``Fantastic!'' cried the Number Devil.``So if you take a few more
steps to end up back where you started, and turn once more with
angle gamma to end up the way you were originally facing, your wheel
would read''---and here the Number Devil drew out his favorite
purple piece of chalk, and began writing in the air\ldots

\begin{eqnarray}
&&(180^\circ - \alpha) + (180^\circ - \beta) + (180^\circ - \gamma)
\nonumber \\
&&\ \ \ \ =\ 540^\circ - (\alpha + \beta + \gamma)
\nonumber \\
&&\ \ \ \ =\ 360^\circ + 180^\circ - (\alpha + \beta + \gamma)
\nonumber \\
&&\ \ \ \ \to\  180^\circ - (\alpha + \beta + \gamma)
\nonumber
\end{eqnarray}

``Fantastic!'' cried Robert. ``Just what we were looking for in
order to find the area of the triangle! I bet if you were on a flat
surface, the wheel would read exactly zero degrees when you
completed the triangle.''

``That's correct---the angles add up to $180^\circ$ then!'' said the
Number Devil. ``But the fun doesn't stop here. One of the beautiful
properties of spherical geometry is that you can calculate the area
of any shape just from the interior angles, as long as its sides are
segments of Great Circles. So if you walked around an $N$-sided
polygon with interior angles $\alpha_1$, $\alpha_2$,\ldots,
$\alpha_N$ you would have an accumulated angle of:''

\begin{equation}
\theta = (N-2) \cdot 180^\circ
- (\alpha_1 + \alpha_2 + \ldots + \alpha_N)
\nonumber
\end{equation}

``Brilliant!'' cried Robert.  ``That's the difference between the
sum of all the polygon's angles, and the sum it would have were it
flat---$180^\circ$ times two less than the number of sides.''

``So,'' the Number Devil continued, ``the area would just be:''

\begin{equation}
A = R^2 \cdot (-\theta) \cdot \left(\frac{2\pi}{360} \right)
\nonumber
\end{equation}

``This remains true for an arbitrary smooth closed curve on the
sphere, which is made of tiny Great Circle segments. Using your
wheel device, you can still read off the accumulated angle after
returning to your starting point, and this formula will still give
you the area enclosed by that curve!'' said the Number Devil.

``Brilliant!'' cried Robert again.

``Even more brilliant is how they used to make these calculations in
the old days, before they could travel to a conveniently sized
asteroid," said the Number Devil. ``In 1851, a scientist named
L{\'e}on Foucault realized that any point on the Earth is rotating
at a high speed around a circle of fixed latitude---a smooth curve
of the type I was just talking about. So, if someone were to stand
still, holding their wheel while the Earth rotated 360 degrees, they
would accumulate an angle\ldots''

\begin{equation}
\theta = - \left(\frac{A}{R^2}\right)
\cdot \left(\frac{360}{2\pi} \right)
\nonumber
\end{equation}

``I see,'' said Robert. ``That's just our area formula flipped around.''

``It gets even more clever,'' said the Number Devil. ``Foucault also
calculated that the area enclosed by the circle of latitude phi on
which the person is sitting would equal:''

\begin{equation}
A = 2\pi R^2[1 - \sin(\phi)]
\nonumber
\end{equation}

``Now where on Earth did that come from?'' cried Robert.

``It's not all that complicated, actually, if you have a little bit
of calcification'' said the Number Devil. ``It comes from the
addition of all of the infinitesimally thin circles, of radius $R$
times the cosine of the latitude, for all of the latitudes between
latitude phi and the North Pole.''

``That makes my head hurt,'' Robert sighed.

``But here is the most beautiful part of all!'' cried the Number
Devil. ``If you combine the two formulas I have written above, you
get one elegant queen of a formula!''

``Alright,'' said Robert. ``Let's see it.''

The Number Devil wrote out, very slowly and reverently:

\begin{equation}
\theta = 360^\circ \sin(\phi)
\nonumber
\end{equation}

``Okay,'' said Robert. ``So if you were in Paris, where phi is
$49^\circ$ degrees, you would measure theta as around $272^\circ$.
And if you were at the North Pole, where phi is $90^\circ$, theta
would be $360^\circ$, as the whole Earth rotates underneath you!''

``You've got it!'' said the Number Devil.``However, we're not quite
done with our friend Foucault. Although we can use a device like our
frictionless parallel transporter in your dreams, Robert, it would
be impossible to actually build one for use on Earth.''

Foucault had everything going against him, didn't he,'' said Robert.
``What did he do instead?''

``He used a very large, heavy pendulum," said the Number Devil. "To
be specific, he hung a 27 kilogram weight from a 67 meter long wire
from the dome of the Pantheon in Paris. This pendulum was able to
swing back and forth, almost frictionlessly, for several hours.''

``And the pendulum must have oscillated!" said Robert. "I think I've
seen one of these at the Science Museum. It doesn't just swing back
and forth, it slowly rotates---you can even use it as a clock!''

``That's right,'' said the Number Devil. ``Except if you were at the
equator\ldots''

``Well then it would just swing back and forth. But that wouldn't be
interesting at all,'' said Robert.

``You've got it, my boy,'' said the Number Devil. ``Now, before I
let you go back to sleep, I want you to take one more look at the
soccer ball that started this lesson in spherical geometry.  Notice
that the pentagons and hexagons covering the soccer ball are all
examples of Great Circle polygons, each covering some part of the
surface of the sphere.''

``Yes,'' said Robert, ``and now I know that the interior angles of
the pentagons add up to more than $540^\circ$, and those of the
hexagons add up to more than $720^\circ$, and the difference is
proportional to their area.  What more is there to know?''

``Well,'' said the Number Devil,''  I want you to think about all of
the polygons together, not one at a time. There is a beautiful
formula, discovered by the great Swiss mathematician Leonhard Euler
around 1740, that relates the number of vertices, faces and edges of
such a polygonal covering of a sphere:''

\begin{equation}
V + F - E = 2
\nonumber
\end{equation}

Robert squinted and circled around the soccer ball, counting up the
various features.  ``Yes, it seems to work.'' he said. ``There are
twelve pentagons and twenty hexagons, and I count sixty vertices,
thirty-two faces and ninety edges, and indeed sixty plus thirty-two
minus ninety is exactly two. Phew! But what does this have to do
with everything else we have talked about?''

``Well,'' said the Number Devil, ''let's think about applying our
angle formula to all of the faces at once.  What is their total
area?''

``The total area is four pi times the square of the radius.'' said
Robert. ``That's just a special case of the calcified formula where
the angle phi is the south pole.''

``Wow!'' said the Number Devil, pulling out his piece of purple
chalk. ``You have been paying attention!  Now let's write that total
area as an angle sum, adding up all of the contributions from all of
the faces:''

\begin{widetext}
\begin{equation}
4\pi = \left(\frac{2\pi}{360} \right) \cdot
\begin{array}{c}
\mbox{\sl SUM OVER ALL FACES} \\
f = 1,2,\ldots,F
\end{array}
\Big\{\alpha_{f,1} + \alpha_{f,2} + \ldots +
\alpha_{f,V_f} - 180 (V_f-2) \Big\}
\nonumber
\end{equation}
\end{widetext}

Robert frowned in concentration.  ``So, the radius-squared has
canceled on both sides, and the number of vertices $V_f$ belonging
to face $f$ is five for the pentagons, and six for the hexagons.''

``Right!'' said the Number Devil. ``But now let's be really tricky,
and notice that this sum is adding up every single angle on the
soccer ball.  So, we can reorganize it, to first add up every angle
surrounding each vertex---instead of each face---and then add up the
result over all vertices---instead of over all faces.  Do you see
what I'm getting at?''

\centerline{\includegraphics[bb = 0 0 378 412,width=3.0in]{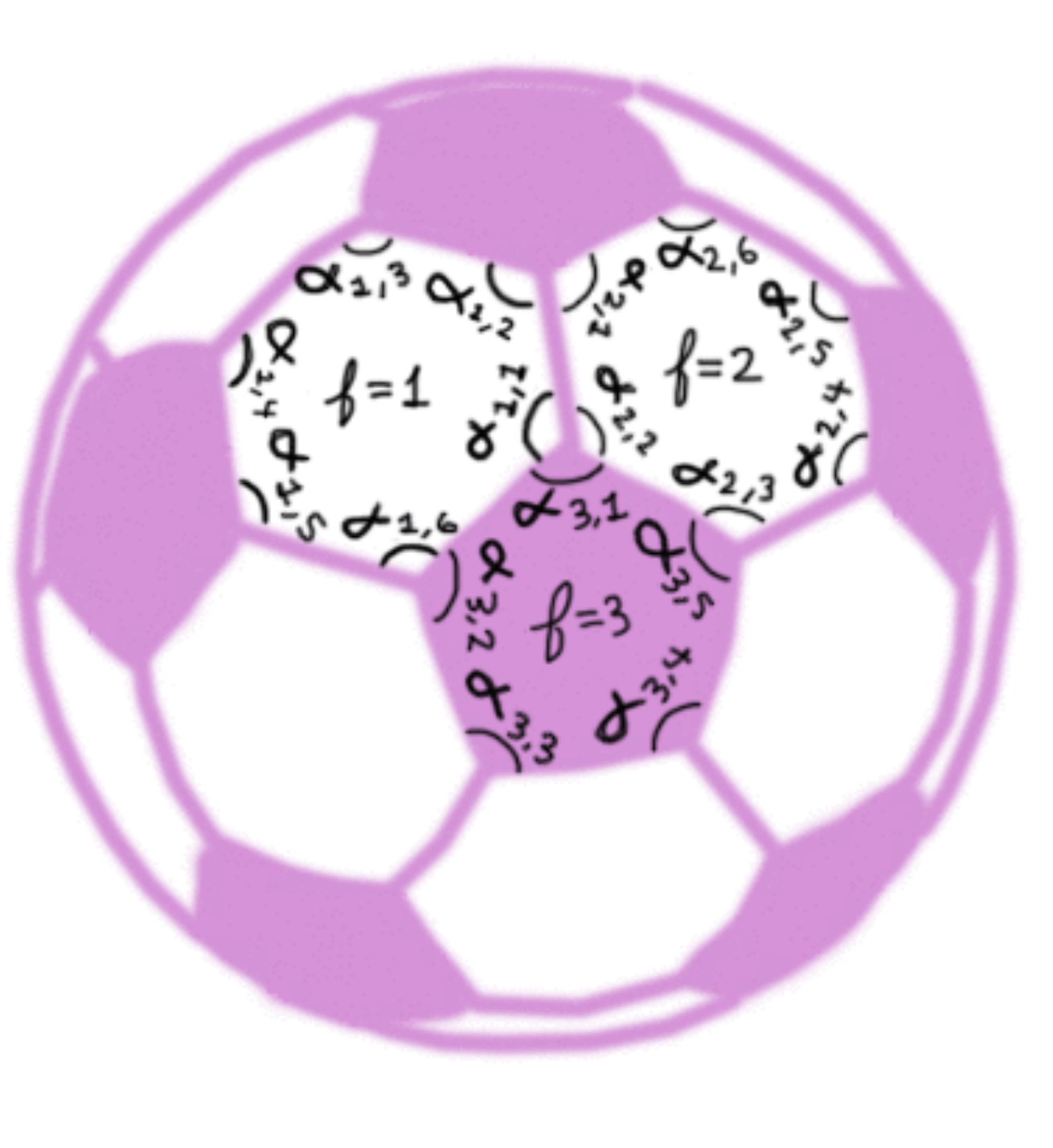}}

Desperate to go back to sleep, Robert concentrated extra hard, and
suddenly saw the pattern.  ``Yes,'' He said, ``I see it now!  Each
has vertex two hexagons and one pentagon around it, and the three
angles must add up to $360^\circ$.  So the sum of all the angles on
the soccer ball is just $360^\circ$ times the total number of
vertices $V$.''

The Number Devil was impressed. ``Wonderful!  So now I can simplify
the relation to:''

\begin{equation}
4\pi = 2\pi V + 2\pi F - \pi (V_1 + V_2 + \ldots + V_F)
\nonumber
\end{equation}

``Well,'' said Robert, ''that's really looking a lot nicer.  But
what do we do with that last sum?  Different faces have different
numbers of vertices, so I don't see how you can simplify it any
further.''

``Ah,'' said the Number Devil in a superior tone, ``but the trick
now is to realize that each polygon has the same number of edges as
vertices: $E_f$ is the same as $V_f$.  Now what's special about
edges?''

``Ummm\ldots'' hesitated Robert. ``Well, an edge is where two faces
meet.''

``Exactly!'' said the Number Devil. ``So when you add up the number
of edges over all the faces, you are counting each edge exactly
twice---once for each of its adjacent faces.  The result of that
last sum is therefore twice the total number of edges, $E$, and we
finally arrive at:''

\begin{equation}
4\pi = 2\pi(V + F - E)
\nonumber
\end{equation}

``Wow!'' exclaimed Robert.  ``If you divide both sides by twice pi
you just get Euler's formula.''

``Yes,'' said the Number Devil, ``and if you think about it, you
will see that none of our steps actually used the fact that we were
dealing with a soccer ball, so this formula works for any polygonal
covering of the sphere. Also, since this formula depends only on the
total number of vertices, faces and edges, we are free to distort
the spherical polygon covering into more standard solid shapes, like
pyramids and cubes, with straight edges and planar faces, without
changing the formula.  In fact, the soccer ball is just an
icosahedron formed from 20 equilateral triangles, with its corners
filed down to make the pentagons.  That way it rolls better!

``You see also that the number `2' relates to something very
special---it comes from the area of the sphere. This is such an
ingenious result that this `2' is called the `genius of the sphere.'
If you were to cover other solids with polygons, you would get
different amounts of genius. For example, if you tried this on a
donut or inflatable inner tube, you would need to replace the `2' by
a `0'---which is logical given that those shapes look like zeros.''

``Ah!'' said Robert. ``Then if I glue two donuts together to make a
figure eight, I bet you should replace the `0' by an `8' for that
type of solid!''

``Er, no.'' said the Number Devil. ``Good guess, but actually you
replace it by a `$-2$'.  We'll have to get a `handle' on that on
another night.

``And now, I think that it really is time to let you go back to
sleep.''

``Finally,'' said Robert. As he drifted off, Robert heard the Number
Devil say, ``And tomorrow we'll learn about hyperbolic geometry and
Klein bottles and Nikolai Ivanovich Lobachevsky\ldots''

Robert awoke the next morning with his head spinning like the Earth,
thinking that if he were to draw all possible lines connecting his
nose, mouth, and eyes, he could distort his head into an tetrahedron
with four faces, four vertices, and six edges\ldots.

\bigskip

\centerline{THE END}

\end{document}